\documentclass[a4paper,11pt]{amsart}
\usepackage{amsmath}
\usepackage{amssymb}
\usepackage{amstext}
\usepackage{amsbsy}
\usepackage{amsopn}
\usepackage{amsthm}
\usepackage{amscd}
\usepackage{amsxtra}
\usepackage{amsfonts}
\usepackage{ifthen}
\usepackage{xspace}
\usepackage{fancyheadings}
\usepackage{layout}
\usepackage{hyperref}
\newcommand{\meinehomepage}{\href{http://www-ifm.math.uni-hannover.de/~bothmer/}
                     {http://www-ifm.math.uni-hannover.de/$\tilde{\,\,}$bothmer}
                           }
\newcommand{\preprintserver}[2]{\href{http://xxx.lanl.gov/abs/math/#2}{#1/#2}}

\input xy
\xyoption{matrix}
\xyoption{arrow}
\xyoption{curve}
\newdir{ (}{{}*!/-5pt/@^{(}}
\newcommand{\xycenter}[1]{\begin{center}
                          \mbox{\xymatrix{#1}}
                          \end{center}
                         }

\newboolean{xlabels}
\newcommand{\xlabel}[1]{
                        \label{#1}
                        \ifthenelse{\boolean{xlabels}}
                                   {\marginpar{#1}}
                                   {}
                       }

\newcommand{\PZ}{\mathbb{P}}

\newcommand{\CC}{\mathbb{C}}

\newcommand{\PP}{\mathbb{P}}
\newcommand{\GG}{\mathbb{G}}

\newcommand{\sE}{{\mathcal E}}

\newcommand{\sO}{{\mathcal O}}

\newcommand{\sQ}{{\mathcal Q}}

\newcommand{\from}{\leftarrow}

\newboolean{probleme}
\newcommand{\problem}[1]
           {\ifthenelse{\boolean{probleme}}
                       {{\bf(PROBLEM: #1)\bf}}
                       {}
           }
\newboolean{zukuenftiges}
\newcommand{\zukunft}[1]
           {\ifthenelse{\boolean{zukuenftiges}}
                       {{\bf(AUSBAUM\"OGLICHKEIT: #1)\bf}}
                       {}
           }
\newboolean{extras}
\newcommand{\extra}[1]
           {\ifthenelse{\boolean{extras}}
                       {{\bf EXTRA #1 EXTRA\bf}}
                       {}
           }

\newboolean{ignore}
\newcommand{\ignore}[1]
           {\ifthenelse{\boolean{ignore}}
                       {{\bf IGNORE #1 IGNORE\bf}}
                       {}
           }

\DeclareMathOperator{\Img}{Im}

\DeclareMathOperator{\rank}{rank}

\DeclareMathOperator{\cliff}{cliff}
\DeclareMathOperator{\Syz}{Syz}
\DeclareMathOperator{\Gensyz}{Gensyz}

\DeclareMathOperator{\GL}{GL}

\theoremstyle{plain}
\begingroup
\newtheorem{thm}{Theorem}%[subsection]
\newtheorem{cor}[thm]{Corollary}
\newtheorem{lem}[thm]{Lemma}

\newtheorem{prop}[thm]{Proposition}
\newtheorem{conj}[thm]{Conjecture}

\newtheorem{question}[thm]{Question}
\numberwithin{thm}{subsection} 
\endgroup

\newtheorem*{thm*}{Theorem}
\newtheorem*{conj*}{Conjecture}
\newtheorem*{verm*}{Vermutung}

\theoremstyle{definition}
\newtheorem{defn}[thm]{Definition}
\newtheorem{rem}[thm]{Remark}
\newtheorem{example}[thm]{Example}

\numberwithin{equation}{section}

\newcommand{\nosubsections}{\renewcommand{\thethm}{\thesection.\arabic{thm}}
                            \setcounter{thm}{0}
                           }

\newcommand{\cref}[3]{(\ref{#1}, #2 \ref{#3})}

\date{\today}

\usepackage{amssymb,amsmath}
\setboolean{probleme}{true}
\setboolean{xlabels}{false}

\newcommand{\secemail}{
\setlength{\unitlength}{1pt}
bothmer
\begin{picture}(0,1)
\put(0,0){m}
\put(-5,0){@}
\end{picture}
ath.uni-hannover.de}

\begin{document}

\title{Generic Syzygy Schemes}

\address{Institiut f\"ur Mathematik\\ 
          Universit\"at Hannover\\ 
          Welfengarten 1\\ 
          D-30167 Hannnover 
         }

\email{\secemail}

\urladdr{\meinehomepage}
\thanks{Supported by the Schwerpunktprogramm ``Global Methods in Complex
        Geometry'' of the Deutsche Forschungs Gemeinschaft}
        
\author{Hans-Christian Graf v. Bothmer}

\begin{abstract}
For a finite dimensional vector space $G$ we define the $k$-th generic syzygy scheme 
$\Gensyz_k(G)$ by explicit equations. If $X \subset \PZ^n$ is cut out by quadrics and $f$ is a $p$-th syzygy of rank $p+k+1$ we show that the syzygy scheme $\Syz(f)$ of $f$ is a cone
over a linear section of $\Gensyz_k(G)$. We also give a geometric description of $\Gensyz_k(G)$ for
$k=0,1,2$, in particular $\Gensyz_2(G)$ is the union of a Pl\"ucker embedded Grassmannian and a linear space. From this we deduce that every smooth, non-degenerate projective curve $C \subset \PZ^{n}$ which is cut out by quadrics and has a $p$-th linear syzygy of rank $p+3$ admits a rank $2$ vector bundle $\sE$ with $\det  \sE = \sO_C(1)$ and
$h^0(\sE) \ge p+4$.
\end{abstract} 

\maketitle

\newcommand{\dual}{^*}

%%%%%%%%%%%%%%%%%%%%%%%%%%%%%%%%%%%%%%%%%%%%
\section{Introduction}
%%%%%%%%%%%%%%%%%%%%%%%%%%%%%%%%%%%%%%%%%%%%
\nosubsections

Let $X \subset \PP^n$ be a projective variety that is cut out by quadrics. One can then look at the linear
strand of its minimal free resolution and ask whether a $p$-th linear syzygy $f$ carries some geometric information about $X$. For this purpose Ehbauer \cite{Ehbauer} introduced the syzygy scheme $\Syz(f)$, which is cut out by the quadrics involved in $f$. The syzygy scheme always contains $X$ and can be explicitly calculated in some cases. Ehbauer studied this construction when $X$ is a set of points in uniform position.  

Another geometric invariant of a $p$-th syzygy $f$ is the space $G\dual$ of linear forms
involved in $f$. Its dimension is called the rank of $f$. Interesting syzygy varieties often arise from syzygies of low rank.  

In \cite{Sch91} Schreyer observed that for $p=1$ the syzygy scheme $\Syz(f)$ is always a cone over a linear section of a generic syzygy scheme $\Gensyz_k$ with $k = \rank f - 2$ and gave explicit equations for $\Gensyz_k$ in this case.  Eusen and Schreyer found a geometric description of these schemes for \(k \in \{0, \dots, 4\}\) and $p=1$ in \cite{EusenSchreyer}. 

In this paper we define more general generic syzygy schemes $\Gensyz_k(G)$ by explicit equations depending on a finite dimensional vector space $G$. With these schemes we prove:

\newtheorem*{thm-cone}{Theorem \ref{t-cone}}

\begin{thm-cone}
Let $I \subset R$ be a homogeneous ideal generated by quadrics and $f$ a $p$-th rank $p+k+1$ linear syzygy of $I$. Then
the syzygy scheme $\Syz(f)$ is isomorphic to a cone over a linear section of $\Gensyz_k(G)$ where $G$ is  the  space of $(p-1)$-st syzygies involved in $f$.
\end{thm-cone}

We also obtain a geometric description of $\Gensyz_k(G)$ for $k=0,1,2$ and arbitrary $G$. 
We show that $\Gensyz_0(G)$ is always
the union of a hypersurface with a point and that $\Gensyz_1(G)$ is a Segre-embedded $\PP^1 \times \PP^{\rank f -1}$. The main new result of this paper is

\newtheorem*{thm-k2}{Theorem \ref{t-k2}}

\begin{thm-k2}
Let $G$ be a $g$ dimensional vector space, then
\[
	\Gensyz_2 (G) =  \GG(\CC \oplus G\dual,2)  \cup \PP\bigl(\bigwedge^2 G\dual\bigr) 
	\subset \PP\bigl(G\dual \oplus \bigwedge^2 G\dual\bigr),
\]
where $\GG(\CC \oplus G\dual,2)$ is the Grassmannian of two dimensional quotient spaces
of $\CC \oplus G\dual$.
Moreover the second generic syzygy ideal $I$ of $G$ is reduced and saturated.
\end{thm-k2}

The geometric descriptions of $\Gensyz_k(G)$ allow us to draw a number of conclusions:

\theoremstyle{plain}
\newtheorem*{cor-reduced}{Corollary \ref{c-red}}

\begin{cor-reduced}
Let $X \subset \PZ^n$ be a projective variety, 
$I_X$ generated by quadrics and $f\in F_p$ a $p$-th syzygy of rank 
$p+1$. Then $X$ is either contained in a hyperplane or reducible.
\end{cor-reduced}

This result seems to be well known, but we include it since it follows directly from our methods.

\newtheorem*{cor-scrolls}{Corollary \ref{c-scrolls}}

\begin{cor-scrolls}
Let $X \subset \PZ^n$ be a non-degenerate irreducible projective variety, 
$I_X$ generated by quadrics and $f\in F_p$ a $p$-th syzygy of rank 
$p+2$. Then the  syzygy scheme $\Syz(f)$ of $f$ is a scroll of degree $p+2$ and
codimension $p+1$.
\end{cor-scrolls}

In particular  a $p$-th syzygy  of rank $p+1$ 
implies the existence of a special pencil $|D|$ on $X$ 
cut out by the fibers of the scroll. If $X$ is a canonical curve $|D|$ has low Clifford index. These pencils are the ones that play a role in Green's conjecture \cite{GreenKoszul}. Our corollary above is therefore probably well known to experts in this field. 

Our main new geometric result is

\newtheorem*{thm-rank2}{Theorem \ref{t-rank2}}

\begin{thm-rank2} 
Let $C \subset \PP^n$ be a smooth, irreducible non-degenerate curve. 
If $C$ is cut out by quadrics and has a $p$-th
syzygy $f$ of rank $p+3$, then there exists a rank $2$ vector bundle $\sE$ on $C$ with
$\det \sE = \sO_C(1)$ and $h^0(\sE) \ge p+4$.
\end{thm-rank2}

In the case of a canonical curve these are rank $2$ budles with canonical determinant.

One can also use the methods of this paper to construct the Mukai-Lazarsfeld bundle on a $K3$ surface directly from a syzygy $f$. This is the vector bundle that played a central role in Voisin's proof of Green's conjecture \cite{VoisinK3}, \cite{VoisinOdd}. The Grassmannian used by Voisin in her proof is dual to the Grassmannian obtained as the generic syzygy scheme of $f$.

This paper is structured as follows. In Section \ref{s-syzygies} we recall what we need about syzygies, syzygy ideals and syzygy schemes. In Section \ref{s-generic} we define the generic syzygy varieties and show that
every syzygy scheme is a cone over a linear section of a generic syzygy scheme.  In the last three sections we describe the $k$-th generic syzygy varieties for $k=0,1,2$ geometrically and study syzygies of rank $p+1$, $p+2$ and $p+3$.

%%%%%%%%%%%%%%%%%%%%%%%%%%%%%%%%%%%%%%%%%%%%
\section{Syzygies, Syzygy Ideals and Syzygy Schemes} \xlabel{s-syzygies}
%%%%%%%%%%%%%%%%%%%%%%%%%%%%%%%%%%%%%%%%%%%%
\nosubsections

For the purpose of this paper let $R = \CC[x_0,\dots,x_n]$ be the homogeneous coordinate ring of $\PP^n$. With $R(-i)$ we denote $R$ with its grading shifted, i.e. $R(-i)_j = R_{j-i}$. Often we abbreviate the space of linear polynomials $R_1 \subset R$ by $V$ and write $\PP^n = \PP(V)$ using the Grothendieck notation.

\begin{defn} \xlabel{d-linearsyz}
Let $I \subset R$ be a homogeneous ideal, generated by quadrics, and
\[
	F_\bullet \colon I \from F_0 \otimes R(-2) \from \dots \from F_r \otimes R(-r-2)
\]
the linear part of the minimal free resolution of $I$. The elements of $F_i$ are called
$i$-th {\sl linear syzygies} of $I$. 
\end{defn}

%\begin{prop}
%In the situation of Definition \ref{d-linsyz} we have
%\[
%	F_p = \ker \Bigl( I_2 \otimes \bigwedge^p V \to I_3 \otimes \bigwedge^{p-1} V \Bigr). 
%\]
%\end{prop}

%\begin{proof}
%Koszulcohomology. \cite{??}
%\end{proof}

\begin{defn} \xlabel{d-rank}
Let $I \subset R$ be a homogeneous ideal, generated by quadrics and $f \in F_p$ a $p$-th 
linear syzygy. We define the {\sl space of $(p-1)$-st linear syzygies involved in $f$}
as the smallest vector space
$G \subset F_{p-1}$ such that there is a commutative diagram
\xycenter{
			F_{p-1} \otimes R(-p-1) & F_p \otimes R(-p-2) \ar[l] \\
			G \otimes R(-p-1) \ar[u] & f \otimes R(-p-2). \ar[l] \ar[u]
		}
We define the {\sl rank} of $f$ as the dimension of $G$.

The above diagram extends to a map from the Koszul complex of $G$ to the linear strand of $I$:
\xycenter{
			I 
			& F_{0} \otimes R(-2) \ar[l] 
			& \dots \ar[l] 
%			& F_{p-1} \otimes R(-p-1) \ar[l] 
			& F_p \otimes R(-p-2) \ar[l] \\
			\bigwedge^{p+1} G \otimes R(-1) \ar[u]
			& \bigwedge^{p} G \otimes R(-2) \ar[u] \ar[l] 
			& \dots \ar[l] 
%			& G \otimes R(-p-1) \ar[u] \ar[l] 
			& f \otimes R(-p-2). \ar[l] \ar[u]\\
		}
The image of $\bigwedge^p G$ in $I$ is called the {\sl syzygy ideal} $I_f$ of $f$.		
\end{defn}

\begin{rem}
Observe that by dualizing and twisting the morphism $$G \otimes R(-p-1) \from f \otimes R(-p-2)$$
from above, $G\dual$ is exhibited
as a space of linear forms on $\PP^n$. We therefore call $G\dual$ the {\sl space of
linear forms involved in $f$}.
\end{rem}

\begin{lem} \xlabel{l-nonzero}
In the map of complexes of Definition \ref{d-rank} all vertical maps are nonzero.
\end{lem}

\begin{proof}
Suppose there exists an integer $k$ such that in the diagram
\xycenter{
			F_{k-1} \otimes R(-k-1)
			& F_{k} \otimes R(-k-2)\ar[l] \\
			\bigwedge^{p-k+1} G \otimes R(-k-1) \ar[u]^{\varphi_{k-1}}
			& \bigwedge^{p-k} G \otimes R(-k-2) \ar[u]^{\varphi_{k}} \ar[l] \\
		}
the morphism $\varphi_{k-1}$ is zero, but $\varphi_{k}$ is nonzero. Then the image
of $\varphi_k$ is a free summand of $F_{k} \otimes R(-k-2)$ which maps to zero in the linear strand of the minimal free resolution of $I$. This contradicts the minimality of the resolution.
\end{proof}

\begin{cor}
Let $f$ be a $p$-th linear syzygy of $I \subset R$. Then $\rank f \ge p+1$. 
\end{cor}

\begin{proof}
If $\rank f \le p$ then $\bigwedge^{p+1} G$ vanishes and the first vertical map of the map of complexes in Definition \ref{d-rank} would have to be zero.
\end{proof}

\begin{defn}
Let $I \subset R$ be an ideal generated by quadrics, $f \in F_p$ a $p$-th linear syzygy and $I_f$
the syzygy ideal of $f$. Then the vanishing set $\Syz(f) = V(I_f)$ is called the syzygy scheme associated
to $f$. 
\end{defn}

\begin{rem} \xlabel{r-notall}
Observe that $\Syz(f) \subset \PP^n$ is always a strict subset, since the syzygy ideal $I_f$ is
never empty by Lemma \ref{l-nonzero}.
\end{rem}

%%%%%%%%%%%%%%%%%%%%%%%%%%%%%%%%%%%%%%%%%%%%
\section{Generic syzygy schemes} \xlabel{s-generic}
%%%%%%%%%%%%%%%%%%%%%%%%%%%%%%%%%%%%%%%%%%%%
\nosubsections

\begin{defn}
Let $G$ be a vector space of dimension $g$ and consider the 
ring $S = \CC[G\dual \oplus \bigwedge^k G\dual]$. The ideal $I$ defined by the natural inclusion
\[
	I = \bigwedge^{k+1} G\dual 
		\subset G\dual   \otimes \bigwedge^k G\dual 
		\subset S^2 \Bigl(G\dual \oplus \bigwedge^k G\dual\Bigr) \subset S
\]
is called the $k$-th generic syzygy ideal of $G$. Its vanishing set $\Gensyz_k(G)$ is called
the $k$-th generic syzygy scheme of $G$.
\end{defn}

\ignore{       %% don't consider cubics at this moment
\begin{prop}
Let $I$ be the $k$-th generic syzygy ideal of $G$. Then the vector space $I_3$ of cubics
in $I$ decomposes in irreducible representations of $\GL(g)$. In this decomposition the
representation $\bigwedge^{k+2}$ does not occur. 
\end{prop}

\begin{proof}
First we have the decomposition
\[
	S_3 = S^3 G \oplus 
		     S^2 G \otimes \bigwedge^k G \oplus 
		     G \otimes S^2 \Bigl( \bigwedge^k G \Bigr)
		     \oplus S^3 \Bigl( \bigwedge^k G \Bigr).
\]
If $k=1$ then the irreducible decomposition of $S_3$ is 
\[
	S_3 = 4 S^3 G \oplus 2S^{2,1} G 
\]
so no summand $\bigwedge^3 G$ occurs. If $k \ge 2$ the four summands of $S_3$
have weights $3 < 2+k < 1 + 2k < 3k$ and $\Lambda^{k+2} G$ can only occur in the second
summand. But this summand decomposes into the irreducible representations $\bigwedge^{k,1,1} G$
and $\bigwedge^{k+1,1} G$. Consequently $\Lambda^{k+2}$ does not occur in the irreducible
decomposition of $S_3$ and therefore also not in the decomposition of $I_3 \subset S_3$.
\end{proof}
} %%

\begin{prop} \xlabel{p-subkoszul}
Let $I$ be the $k$-th generic syzygy ideal of $G$. Then the linear strand of $I$ has
the last $g-k$ steps of the Koszul complex associated to $G\dual$ as a natural subcomplex, i.e
we have a commutative diagram:
\xycenter{
			I 
			& F_{0} \otimes S(-2) \ar[l] 
			& \dots \ar[l] 
%			& F_{p-1} \otimes R(-p-1) \ar[l] 
			& F_{g-k-1} \otimes S(-g+k-1) \ar[l] \\
			%\bigwedge^{k} G\dual \otimes R(-1) \ar@{ (->}[u]
			& \bigwedge^{k+1} G\dual \otimes S(-2) \ar@{ (->}[u] %\ar[l] 
			& \dots \ar[l] 
%			& G \otimes R(-p-1) \ar[u] \ar[l] 
			& \bigwedge^g G\dual \otimes S(-g+k-1) \ar[l] \ar@{ (->}[u]\\
		}
\end{prop}

\begin{proof}
The inclusion 
$	I = \bigwedge^{k+1} G\dual 
		\subset G\dual   \otimes \bigwedge^k G\dual 
		\subset S^2 \Bigl(G\dual \oplus \bigwedge^k G\dual\Bigr) \subset S$
induces a commutative diagram of free $S$-modules
\xycenter{
			I
			& F_{0} \otimes S(-2) \ar[l] \\
			\bigwedge^{k} G\dual \otimes S(-1) \ar[u]
			& \bigwedge^{k+1} G\dual \otimes S(-2). \ar@{=}[u] \ar[l]\\
		}
The top arrow is resolved by the minimal free resolution of $I$ and the bottom arrow
by the rest of the Koszul complex. Since both complexes are exact and minimal, the maps above lift to a map of complexes. This map is injective in each new step since it is injective in the $F_0$ step.
For degree reasons, the image of this map of complexes must lie in the linear strand of $I$. 
\end{proof}

\begin{cor}
The $k$-th generic syzygy scheme of $G$ has a natural $1$-dimensional space of
rank $g$ linear syzygies in step $g-k-1$. The space of $(g-k-2)$-nd syzygies involved in 
anyone of these is isomorphic to $G$.
\end{cor}

\begin{proof}
The $(g-k-1)$-st syzygies given by Proposition \ref{p-subkoszul} have rank at most $g$ since
$\bigwedge^{g-1} G\dual \cong G$ has dimension $g$. The rank of these syzygies cannot be smaller,
since the last map of the Koszul complex is surjective in degree $g-k$.
\end{proof}

\begin{thm} \xlabel{t-cone}
Let $I \subset R$ be a homogeneous ideal generated by quadrics and $f$ a $p$-th rank $p+k+1$ linear syzygy of $I$. Then
the syzygy scheme $\Syz(f)$ is isomorphic to a cone over a linear section of $\Gensyz_k(G)$ where $G$ is  the  space of $(p-1)$-st syzygies involved in $f$.
\end{thm}

\begin{proof}
We have
the map of complexes
\xycenter{
			R
			& F_{0} \otimes R(-2) \ar[l] 
			& \dots \ar[l] 
%			& F_{p-1} \otimes R(-p-1) \ar[l] 
			& F_p \otimes R(-p-2) \ar[l] \\
			\bigwedge^{p+1} G \otimes R(-1) \ar[u]^{\alpha}
			& \bigwedge^{p} G \otimes R(-2) \ar[u] \ar[l] 
			& \dots \ar[l] 
%			& G \otimes R(-p-1) \ar[u] \ar[l] 
			& f \otimes R(-p-2) \ar[l] \ar[u]\\
		}
from Definition \ref{d-rank}. Consider the map
\[
	\varphi \colon G\dual \oplus \bigwedge^k G\dual \to V
\]
given by mapping the elements of $G\dual$ to their corresponding linear forms and
the elements of $\bigwedge^k G \dual = \bigwedge^{p+1} G $ to their images under the map $\alpha$. The induced diagram
\xycenter{
			R
			& F_{0} \otimes R(-2) \ar[l] \\
			\bigwedge^{p+1} G \otimes R(-1) \ar[u]^{\alpha}
			& \bigwedge^{p} G \otimes R(-2) \ar[u] |!{[r];[lu]}\hole \ar[l]  
			& S \ar[ull]
			& \bigwedge^{k+1} G\dual \otimes S(-2) \ar[l] \ar[ull]\\
			&
			&\bigwedge^{k} G\dual \otimes S(-1) \ar[u] \ar[ull]
			& \bigwedge^{k+1} G\dual \otimes S(-2) \ar@{=}[u] \ar[l] \ar[ull] |!{[l];[lu]}\hole \\
		}
and its degree $2$ part
\xycenter{
			S^2 V
			& F_{0}  \ar[l] \\
			\bigwedge^{p+1} G \otimes V \ar[u]^{\alpha}
			& \bigwedge^{p} G  \ar[u] |!{[r];[lu]}\hole \ar[l]  
			& S^2(G\dual \oplus \bigwedge^k G\dual) \ar[ull]
			& \bigwedge^{k+1} G\dual \ar[l] \ar[ull]\\
			&
			&\bigwedge^{k} G\dual \otimes (G\dual \oplus \bigwedge^k G\dual) \ar[u] \ar[ull]
			& \bigwedge^{k+1} G\dual \ar@{=}[u] \ar[l] \ar@{=} [ull] |!{[l];[lu]}\hole \\
		}
shows
that $\varphi$ maps the $k$-th generic syzygy ideal surjectively to the
syzygy ideal $I_f$ of $f$. 
Projectively the image of $\varphi$ defines a linear subspace 
$$\PP(\Img \varphi) \subset \PZ(G\dual \oplus \bigwedge^k G\dual).$$
The calculation above shows
that $\Syz(f)$ is a cone over $\PP(\Img \varphi) \cap \Gensyz_k(G)$ with vertex
$V(\Img \varphi) \subset \PP(V)$.
 \end{proof}   
 
 %%%%%%%%%%%%%%%%%%%%%%%%%%%%%%%%%%%%%%%%%%%%
 \section{Reducible Syzygies}
 %%%%%%%%%%%%%%%%%%%%%%%%%%%%%%%%%%%%%%%%%%%%
 \nosubsections
 
 \begin{prop} \xlabel{p-k0}
 Let $G$ be a $g$ dimensional vector space, then
 \[
              \Gensyz_0 (G) \cong \PP(G\dual) \cup \PP(\CC) 
              				\subset \PP(G\dual \oplus \CC),
\]
i.e $\Gensyz_0(G)$ is the union of a hyperplane and a point. Moreover the generic syzygy ideal
of $I$ of $\Gensyz_0(G)$ is reduced and saturated.
\end{prop}

\begin{proof}
The ideal of the hyperplane $\PP(G\dual) \cong\PP^{g-1}$ is generated by the linear forms in $\bigwedge^0 G\dual \cong \CC$. The ideal of the point $\PP(\bigwedge^0 G\dual) \cong \PP^0$ is generated
by the linear forms in $G\dual$. Since the two ideals involve different sets of variables, their
intersection is the same as their product:
\begin{align*}
	I_{\PP^{g-1}} \cap I_{\PP^0} 
	= (G\dual) \cap \bigl(\bigwedge^0 G\dual \bigr) 
	= (G\dual) \cdot \bigl(\bigwedge^0 G\dual \bigr) 
	= (G\dual \otimes \bigwedge^0 G\dual \bigr) 
	= (\bigwedge^1 G\dual)
\end{align*}
This is the $0$-th generic syzygy ideal of $G$.
\end{proof}

\begin{cor} \xlabel{c-red}
Let $X \subset \PZ^n$ be a projective variety, 
$I_X$ generated by quadrics and $f\in F_p$ a $p$-th syzygy of rank 
$p+1$. Then $X$ is either contained in a hyperplane or reducible.
\end{cor}

\begin{proof}
By Theorem \ref{t-cone} and Proposition \ref{p-k0} $\Syz(f)$ is a cone over a linear section of
a hyperplane and a point. Since $\Syz(f)$ can not contain all of $\PP^n$ by Remark \ref{r-notall},
 $\Syz(f) \subset \PP^n$ 
must be the union of a hyperplane and possibly a second linear subspace. Since $X$ is contained
in $\Syz(f)$ it must be either reducible or contained in one of the two linear subspaces.
\end{proof}

\begin{defn}
Let $X \subset \PP^n$ be a projective scheme, whose ideal is cut out by quadrics. 
A $p$-th linear syzygy of $X$ is called {\sl reducible}, if it has rank $p+1$. 
\end{defn}

%%%%%%%%%%%%%%%%%%%%%%%%%%%%%%%%%%%%%%%%%%%%%
\section{Scrollar Syzygies}
%%%%%%%%%%%%%%%%%%%%%%%%%%%%%%%%%%%%%%%%%%%%
\nosubsections

\begin{thm}
Let $G$ be a $g$ dimensional vector space, then
\[
	\Gensyz_1 (G) =  \PP(G\dual) \times \PP^1 \subset  \PP(G\dual \oplus G\dual).
\]
Moreover the second generic syzygy ideal $I$ of $G$ is reduced and saturated.
\end{thm}

\begin{proof}
Observe that $G\dual \otimes (\CC \oplus \CC) = G\dual \oplus G\dual$. We can therefore
consider the Segre embedding
\[
	\PP^{g-1} \times \PP^1 = \PP(G\dual) \times \PP(\CC \oplus \CC) \subset \PP(G\dual \oplus G\dual).
\]
The ideal of $\PP^{g-1} \times \PP^1$ is generated by the Segre quadrics:
\[
	I_{\PP^{g-1} \times \PP^1}
	=  \bigl(\bigwedge^2 G\dual \otimes \bigwedge^2 (\CC \oplus \CC) \bigr)
	=  \bigl(\bigwedge^2 G\dual \bigr)
\]
This is the first generic syzygy ideal of $G$.
\end{proof}

\begin{cor} \xlabel{c-scrolls}
Let $X \subset \PZ^n$ be a non degenerate irreducible projective variety, 
$I_X$ generated by quadrics and $f\in F_p$ a $p$-th syzygy of rank 
$p+2$. Then the  syzygy scheme $\Syz(f)$ of $f$ is a scroll of degree $p+2$ and
codimension $p+1$.
\end{cor}

\begin{proof}
Let $G$ be the $g=p+2$ dimensional space of $(p-1)$-st syzygies involved in $f$. 
By theorem \ref{t-cone} the syzygy scheme $\Syz(f)$ is a linear section
of a cone over $\PZ^{p+1} \times \PZ^1$. 
Since $\PZ^{p+1} \times \PZ^1$ has codimension $p+1$ and degree
$p+2$ in $\PZ(G\dual \oplus \Lambda^1 G\dual)$ we only have to prove that
this intersection is of expected codimension. By Eisenbud 
\cite[Ex. A2.19]{Ei95}
this is the case if the matrix $M$
whose $2 \times 2$-minors cut out $\PZ^{p+1} \times \PZ^1$ remains
$1$-generic after we apply the map
\[
	\varphi \colon G\dual \oplus \Lambda^1 G\dual \to V
\]
from the proof of Theorem \ref{t-cone}. 

If $\varphi(M)$ is not $1$-generic, we can choose bases of $G\dual$ and $\CC \oplus \CC$ such that
$\varphi(M)$ has the form
\[
	M = \begin{pmatrix}
		l_1 & \dots & l_i & l_{i+1} & \dots & l_{g} \\
		a_1 & \dots &a_i & 0 & \dots & 0 \\
	       \end{pmatrix}
\]
with $l_1,\dots, l_{p+1}$  a basis of  $G\dual$  and $a_1,\dots,a_i$ linearly independent.
Since the syzygy ideal $I_f$ cannot be empty by Lemma \ref{l-nonzero}, $i$ has to be at least $1$.
In this situation $I_f$ contains the $2\times 2$ minor
\[
      \det \begin{pmatrix}
           l_1& l_g \\
           a_1 &  0 \\
           \end{pmatrix}
          = l_g \cdot a_1
\]
which implies that $X$ must be reducible or degenerate. This contradicts our assumptions.
\end{proof}

\begin{defn}
Let $X \subset \PP^n$ be a projective scheme, whose ideal is cut out by quadrics. 
A $p$-th linear syzygy of $X$ is called {\sl scrollar}, if it has rank $p+2$. 
\end{defn}

\begin{example} Let $C \subset \PZ^{g-1}$ be a non hyperelliptic canonical curve
of genus $g$ and $|D|$ a pencil of Clifford
index $\cliff(D) = g-p-3$. 
The $p$-th syzygy of $C$ constructed by the method of  Green and Lazarsfeld in \cite{GL1}
is scrollar. 
\end{example}

With the above geometric description of scrollar syzygy varieties one can prove the following
well known converse of the Green-Lazarsfeld construction:

\begin{prop}
\xlabel{clifford}
Let $C \subset \PZ^{g-1}$ be a non hyperelliptic canonical curve
of genus $g$ and $f \in F_p$ a $p$-th scrollar syzygy. Then there exists a 
linear system  $|D|$ on $C$ with Clifford index $\cliff(D) \le g-p-3$.
\end{prop}

\begin{proof}
Let $G\dual$ be the $p+2$ dimensional space of linear forms involved in $f$. Then the syzygy scheme $\Syz(f)$ of
$f$ is a scroll that contains $C$ and has the vanishing set $V(G\dual)$ as a fiber. Set 
$D = C \cap V(G\dual)$. Since $C \subset \PP^{g-1}$ is non degenerate, $D$ is a divisor on $C$. We consider the linear system $|D|$. Since $D$ is cut out by the ruling of $\Syz(f)$ we have $h^0(D) \ge 2$.
Also $h^0(K-D) \ge p+2$ since the linear forms in $G\dual$ cut out canonical divisors of $C$ that contain $D$. Riemann-Roch now gives:
\begin{multline*}
	\cliff D := d - 2r = (h^0(D) - h^0(K-D)-1+g) -2h^0(D) + 2= \\
				           = g + 1 -h^0(D) - h^0(K-D) \ge g + 1 -2 -(p+2) = g-p-3.
\end{multline*}
\end{proof}

\begin{rem}
For general $k$-gonal canonical curves $C$ Green's conjecture is equivalent to the claim that every
step of  the linear strand $C$ contains at least one scrollar syzygy. This was recently shown by Voisin \cite{VoisinK3}, \cite{VoisinOdd}.
\end{rem}

More generally one can make the following conjecture

\begin{conj}[Generic Geometric Syzygy Conjecture]
Let $C \subset \PZ^{g-1}$ be a general canonical curve of genus
$g$. Then for every $p$ the space of $p$-th linear syzygies 
of $C$ is {\sl spanned} by scrollar syzygies.
\end{conj}

This conjecture is known for for $p=1$ when $g \not= 8$ and for $p=2$ when $g=8$
by \cite{HC} and \cite{HCandRanestad}.

%%%%%%%%%%%%%%%%%%%%%%%%%%%%%%%%%%%%%%%%%%%%%
\section{grassmannian syzygies}
%%%%%%%%%%%%%%%%%%%%%%%%%%%%%%%%%%%%%%%%%%%%
\nosubsections

\begin{thm} \xlabel{t-k2}
Let $G$ be a $g$ dimensional vector space, then
\[
	\Gensyz_2 (G) =  \GG(\CC \oplus G\dual,2)  \cup \PP\bigl(\bigwedge^2 G\dual\bigr) 
	\subset \PP\bigl(G\dual \oplus \bigwedge^2 G\dual\bigr),
\]
where $\GG(\CC \oplus G\dual,2)$ is the Grassmannian of two dimensional quotient spaces
of $\CC \oplus G\dual$.
Moreover the second generic syzygy ideal $I$ of $G$ is reduced and saturated.
\end{thm}

\begin{proof}
Observe that $\bigwedge^2 (\CC \oplus G\dual) = G\dual \oplus \bigwedge^2 G\dual$. We can 
therefore consider the Pl\"ucker embedding
\[
	\GG := \GG(\CC \oplus G\dual,2)  \subset \PP\bigl(G\dual \oplus \bigwedge^2 G\dual\bigr)
\]
and the ideal of the Grassmannian $\GG$ which is generated by $4 \times 4$-pfaffians of a skew symmetric matrix. More
precisely:
\[
	I_\GG = \bigl(\bigwedge^4 (\CC \oplus G\dual) \bigr) 
		  = \bigl(\bigwedge^3 G\dual \oplus  \bigwedge^4 G\dual\bigr) 
		  \subset S^2\bigl(G\dual \oplus \bigwedge^2 G\dual\bigr).
\]
On the other hand $\PP(\bigwedge^2 G\dual) \cong \PP^{{g \choose 2}-1} =: \PP$ is cut out by the linear forms
in $G\dual$, so $I_{\PP} = (G\dual)$. To prove the theorem we calculate the
intersection of these two irreducible ideals:
\begin{align*}
	I_{\PP} \cap I_\GG 
	&= (G\dual) \cap \bigl(\bigwedge^3 G\dual \oplus  \bigwedge^4 G\dual\bigr) \\
	&= \bigl((G\dual) \cap \bigl(\bigwedge^3 G\dual\bigr)\bigr) 
		+ \bigl((G\dual) \cap \bigl(\bigwedge^4 G\dual\bigr)\bigr).
\end{align*}
Now the quadrics in the ideal $(G\dual)$ are given by the image of 
$$G\dual \otimes \bigl(G\dual \oplus \bigwedge^2 G\dual\bigr) 
\to S^2\bigl(G\dual \oplus \bigwedge^2 G\dual\bigr),$$ 
i.e
\[
	(I_{\PP} )_2 = S^2 G\dual \oplus G\dual \otimes \bigwedge^2 G\dual
				    = S^2 G\dual \oplus \bigwedge^3 G\dual \oplus \bigwedge^{2,1}G\dual.
\]
This shows that $(\bigwedge^3 G\dual)$ is contained in $(G\dual)$. For the second intersection of ideals notice that $\bigwedge^4 G\dual$ is contained in $S^2(\bigwedge^2 G\dual)$. So the
generators of 
$(\bigwedge^4 G\dual)$ and  $(G\dual)$ involve different sets of variables and the intersection of
the two ideals is
the same as their product:
\begin{multline*}
	(G\dual) \cap \bigl(\bigwedge^4 G\dual\bigr)  
		= (G\dual) \cdot \bigl(\bigwedge^4 G\dual\bigr)
		= \bigl(G\dual \otimes \bigwedge^4 G\dual\bigr) =\\ 
		= \bigl(\bigwedge^5 G\dual \oplus \bigwedge^{4,1} G\dual\bigr)
		\subset G\dual \otimes S^2 \bigl(\bigwedge^2 G\dual\bigr)
		\subset S^3\bigl(G\dual \oplus \bigwedge^2 G\dual\bigr).
\end{multline*}
On the other hand the cubics of $(\bigwedge^3 G\dual)$ contain
$$\bigwedge^3 G\dual \otimes \bigwedge^2 G\dual 
= \bigwedge^5 G\dual \oplus \bigwedge^{4,1} G\dual \subset G\dual \otimes S^2 \bigl(\bigwedge^2 G\dual\bigr).$$
Since these representations occur only once in $G\dual \otimes S^2 (\bigwedge^2 G\dual)$ they must
be the ones that generate the product of ideals above. In total we have shown 
\[
	I_{\PP} \cap I_\GG = \bigl(\bigwedge^3 G\dual\bigr)
\]
which is the second generic syzygy ideal of $G$.
\end{proof}

\begin{defn}
Let $X \subset \PP^n$ be a projective scheme, whose ideal is cut out by quadrics. 
A $p$-th linear syzygy of $X$ is called {\sl grassmannian}, if it has rank $p+3$. 
\end{defn}

\begin{example} \xlabel{e-grass}
Let $X$ be a $K3$-surface of sectional genus $g$ in $\PP^g$ with Picard group generated by a general hyperplane section $H$. Then $X$ has grassmannian $p$-th syzygies for $p \le \frac{g-4}{2}$.
\end{example}

\begin{proof}
$X$ is cut out by quadrics. Since $X$ is irreducible and non-degenerate, $X$ has no reducible syzygies and does not lie on quadrics of rank $2$ or $1$. $X$ can also not lie on a quadrics of rank $4$ or $3$, since  in this case the
rulings of the quadrics would cut out divisors of degree smaller than $H$ on $X$. Hence, because scrolls are cut out by
$2\times 2$ minors of rank at most $4$, $X$ can have no scrollar syzygies.

Now intersect $X$ with a general hyperplane $H$. Then $X \cap H = C \subset \PP^{g-1}$
is a canonical curve whose minimal free resolution is the restriction of the minimal free resolution
of $X$ to $H$. By the construction of Green and Lazarsfeld $C$ has scrollar $p$-th syzygies for 
$p \le \frac{g-4}{2}$. The rank of a syzygy $f$ can fall by at most one when restricting to a general hyperplane (i.e. when the linear form defining $H$ is involved in $f$). Since $X$ has no scrollar
syzygies, the scrollar syzygies of $C$ must come from grassmannian syzygies of $X$.
\end{proof}

We now describe some geometric consequences of grassmannian syzygies. For this
let  $\sQ$ be the universal rank $2$ quotient bundle on the Grassmannian $\GG = \GG(\CC \oplus G\dual,2)$. The global sections of $\sQ$ are
given by $H^0(\GG,\sQ) = \CC \oplus G\dual$.

\begin{lem} \xlabel{l-sectionideal}
Let $s \in H^0(\GG,\sQ)$ be a global section and $I_s$ the ideal of its vanishing locus on $\GG$. Then
$I_s$ is generated by hyperplane sections of $\GG$, more precisely
\[
	I_s = \bigl(s \wedge H^0(\GG,\sQ)\bigr).
\]
\end{lem}

\begin{proof}
Consider the Koszul complex associated to $s$:
\[
	0 \to \sO_\GG \xrightarrow{s} \sQ \to I_s \otimes \bigwedge^2  \sQ \to 0
\]
Taking cohomology shows $\bigl(s \wedge H^0(\GG,\sQ)\bigr) \subset I_s$. Since $\sQ$ is
globally generated, the converse also follows.
\end{proof}

\begin{rem}
Observe that for a section $s \in \CC \subset \CC \oplus G\dual = H^0(\GG,\sQ)$ we have 
$I_\GG + I_{\PP} = I_s$. In other words a grassmannian syzygy $f$ defines up to
a constant a section of $\sQ$.
\end{rem}

\begin{lem} \xlabel{l-injective}
Let $X \subset \PP(V)$ be a projective variety cut out by quadrics, $f$ a $p$-th grassmannian syzygy of $X$, $G$ the space of $(p-1)$-st syzygies involved in $f$, and
$\varphi \colon G\dual \oplus \bigwedge^2 G\dual \to V$ the induced map. Then  the natural
map
\[
	H^0(\GG,\sQ) \to H^0(\GG \cap \PP(\Img \varphi),\sQ|_{\GG \cap \PP(\Img \varphi)} )
\]
is injective. 
\end{lem}

\begin{proof}
By construction $\Img \varphi$ contains $G\dual$ so the non-zero elements
of $G\dual$ are not contained in $I_{\PP(\Img \varphi)}$. On the other hand the vanishing
ideal
\[
	I_s = \bigl(s \wedge (\CC \oplus G\dual)\bigr)
\]
contains the whole space $\CC \wedge G\dual = G\dual $ if $s \in \CC$, or a non-zero element of $G\dual \wedge \CC=G\dual$ if $s \in G\dual$. So $I_s$ can never be contained in $I_{\PP(\Img \varphi)}$
and $H^0(\sQ \otimes I_{\GG \cap \PP(\Img \varphi)/\GG}) = 0$. The proposition then follows
 from the exact sequence
\[
	0 \to \sQ \otimes I_{\GG \cap \PP(\Img \varphi)/\GG}
	   \to \sQ
	   \to \sQ|_{\GG \cap \PP(\Img \varphi)/\GG}
	   \to 0.
\]
\end{proof}

\begin{thm} \xlabel{t-rank2}
Let $C \subset \PP^n$ be a smooth, irreducible non-degenerate curve. 
If $C$ is cut out by quadrics and has a $p$-th
grassmannian syzygy $f$, then there exists a rank $2$ vector bundle $\sE$ on $C$ with
$\det \sE = \sO_C(1)$ and $h^0(\sE) \ge p+4$.
\end{thm}

\newcommand{\opencone}{{Y^\circ}}

\begin{proof}
Let $\Syz(f)$ be the syzygy scheme of $f$. By Theorems \ref{t-cone} and \ref{t-k2} $\Syz(f)$ is a cone over
a linear section of $\GG(p+4,2) \cup \PP^{{p+3 \choose 2} - 1}$. Now $\Syz(f)$ contains $C$ and $C$ is irreducible and non-degenerate, so $C$ must be contained in a cone $Y$ over a linear section of $\GG$. The universal quotient bundle $\sQ$ on $\GG$ restricts to $\GG \cap \PP(\Img \phi)$ and pulls back to a rank $2$ vector bundle $\sQ_\opencone$ on $\opencone = Y\backslash V(\Img \varphi)$.
If $C$ does not intersect the vertex $V(\Img \varphi)$ of $Y$ the restriction of $\sQ_\opencone$ to $C$  is a vector bundle $\sE$.

If $C$ intersects the vertex of $Y$ in a divisor, we consider the blowup $\tilde{Y}$ of $Y$ in the
vertex. $\sQ$ then pulls back to a rank $2$ vector bundle $\sQ_{\tilde{Y}}$ on $\tilde{Y}$. Since $C$ is smooth the strict transform $\tilde{C}$ of $C$ is isomorphic to $C$ and $\sQ_{\tilde{Y}}$ restricts
to a rank $2$ vector bundle $\sE$ on $\tilde{C}\cong C$.

Finally $C$ can not be contained in the vertex of $Y$ since $C$ is non-degenerate. 

By Lemma \ref{l-injective} we  have $h^0(\sQ|_\GG \cap \PP(\Img \phi)) \ge p+4$. These sections extend to $\opencone$. By Lemma \ref{l-sectionideal} the zero loci of sections of $\sQ$ are cut out by linear forms and their closures contain the vertex of $Y$. Since
$X$ is non-degenerate it can not lie in one of these zero loci, so all sections of $\sQ$ descend
to sections of $\sE$.
\end{proof}

\begin{example}
Our method can is some cases also be used to obtain vector bundles on varieties of higher dimension.
Let for example $X \subset \PP^g$ be a $K3$ surface of even sectional genus $g=2k$ whose Picard group is generated by
a general hyperplane section. Then $X$ has a grassmannian $(k-2)$-nd syzygy by
the argument of Example \ref{e-grass}. One can show that in this case the map 
$$\varphi \colon G\dual \oplus \bigwedge^2 G\dual \to V$$ is surjective. Therfore $\Syz(f)$ is not a cone, and $\sQ$ restricts to a rank $2$ vector bundle $\sE$ on $X$ with $\det \sE = \sO_X(1)$ and $h^0(\sE) \ge k+2$. This is the Mukai-Lazarsfeld bundle 
used by Voisin in her proof of Green's conjecture \cite{VoisinK3}.
\end{example}
 
 This example leads us to ask
 
\begin{question}
Let $X \subset \PP^n$ be a surface cut out by quadrics whose Picard group is generated by
a general hyperplane section. Does every step of the linear strand of $X$ contain a grassmannian syzygy?
\end{question}

\begin{rem}
Voisin's Theorem about the syzygies of K3 surfaces in \cite{VoisinK3} %and \cite{VoisinOdd} 
prove that the answer to this question is
"yes" in the case of K3 surfaces $X \subset \PP^g$ with sectional genus $g=2k$.
\end{rem}

Even more generally we ask

\begin{question}
Let $X \subset \PP^n$ be a surface cut out by quadrics whose Picard group is generated by
a general hyperplane section. Is the space of $p$-th linear syzygies of $X$ {\sl spanned} by grassmannian syzygies?
\end{question}

\begin{rem}
The answer to this question is "yes" for general $K3$ surfaces $X \subset \PP^g$ with sectional genus $g \le 8$ by the methods of
\cite{HCandRanestad}
\end{rem}

%\bibliographystyle{alpha} 
%\bibliography{../bib}                              

\begin{thebibliography}{Ehb94}

\bibitem[Ehb94]{Ehbauer}
S.~Ehbauer.
\newblock Syzygies of points in projective space and applications.
\newblock In F.~Orecchia, editor, {\em Zero-dimensional schemes. Proceedings of
  the international conference held in Ravello, Italy, June 8-13, 1992}, pages
  145--170, Berlin, 1994. de Gruyter.

\bibitem[Eis95]{Ei95}
D.~Eisenbud.
\newblock {\em Commutative Algebra with a View Toward Algebraic Geometry}.
\newblock Graduate Texts in Mathematics 150. Springer, 1995.

\bibitem[ES94]{EusenSchreyer}
F.~Eusen and F.O. Schreyer.
\newblock {A remark to a conjecture of Paranjape and Ramanan}.
\newblock {\href{http://www.math.uni-sb.de/~ag-schreyer/DE/publikationen.html}
  {http://www.math.uni-sb.de/\~{ }ag-schreyer/DE/publikationen.html}}, 1994.

\bibitem[GL84]{GL1}
M.~Green and R.~Lazarsfeld.
\newblock The non-vanishing of certain {Koszul} cohomology groups.
\newblock {\em J. Diff. Geom.}, 19:168--170, 1984.

\bibitem[Gre84]{GreenKoszul}
M.L. Green.
\newblock Koszul cohomology and the geometry of projective varieties.
\newblock {\em J. Differential Geometry}, 19:125--171, 1984.

\bibitem[Sch91]{Sch91}
F.O. Schreyer.
\newblock A standard basis approach to syzygies of canonical curves.
\newblock {\em J. reine angew. Math.}, 421:83--123, 1991.

\bibitem[vB00]{HC}
H.-Chr.~Graf v.~Bothmer.
\newblock {\em Geometrische Syzygien von kanonischen Kurven}.
\newblock Dissertation, Universit\"at Bayreuth, 2000.

\bibitem[vB02]{HCandRanestad}
H.-Chr.~Graf v.~Bothmer.
\newblock Geometric syzygies of {M}ukai varieties and general canonical curves
  with genus $\le 8$.
\newblock \preprintserver{math.AG}{0202133}, 2002.

\bibitem[Voi02]{VoisinK3}
C.~Voisin.
\newblock Green's generic syzygy conjecture for curves of even genus lying on a
  {$K3$} surface.
\newblock {\em J. Eur. Math. Soc. (JEMS)}, 4(4):363--404, 2002.

\bibitem[Voi03]{VoisinOdd}
C.~Voisin.
\newblock Green's canonical syzygy conjecture for generic curves of odd genus.
\newblock \preprintserver{math.AG}{0301359}, 2003.

\end{thebibliography}

\end{document}